\documentclass{amsart}%
\usepackage{amsmath}
\usepackage{graphicx}
\usepackage{amsfonts}
\usepackage{amssymb}
\usepackage{listings}
\usepackage{caption}%
\providecommand{\U}[1]{\protect\rule{.1in}{.1in}}
\makeatletter
\providecommand{\bigsqcap}{\mathop{\mathpalette\@updown\bigsqcup}}
\newcommand*{\@updown}
\makeatother

\theoremstyle{plain}

\newtheorem{algorithm}{Algorithm}

\newtheorem{definition}{Definition}
\newtheorem{example}{Example}

\newtheorem{remark}{Remark}

\numberwithin{equation}{section}
\begin{document}
	\title{Computing the Dimension of a Bipartition Matrix}
	\author{Dawson Freeman}
	\address{Department of Mathematics\\
		Lehigh University\\
		Bethlehem, PA 18015}
	\email{djf321@lehigh.edu}
	\author{Ronald Umble}
	\address{Department of Mathematics\\
		Millersville University of Pennsylvania\\
		Millersville, PA. 17551}
	\email{ron.umble@millersville.edu}
	\date{January 15, 2025}
	\subjclass[2020]{03E05, 05A18, 52B05, 52B11 }
	\keywords{Permutahedron, associahedron, bipartition, bipartition matrix}

\begin{abstract}
	The dimension of a bipartition matrix (BPM) is the sum of the dimensions of
	its indecomposable factors. The dimension of an indecomposable BPM is the
	sum of its row, column, and entry dimensions. To compute these dimensions, we 
	apply four routines of independent interest: (1) Factor a
	bipartition as a product of indecomposables; (2) recover a bipartition from 
	its indecomposable factorization; (3) factor a BPM as a product of 
	indecomposables; and (4) compute the ``transpose-rotation''
	(the column dimension of a BPM is the row dimension of its transpose-rotation).
\end{abstract}

\maketitle

\section{Introduction}

The \emph{permutahedron} $P_{n}$ is an $\left(  n-1\right)  $-dimensional
contractible polytope whose vertices are identified with the $n!$ permutations
of $\underline{n}=\{1,2,\dots,n\}$. In low dimensions, $P_{1}$ is a point,
$P_{2}$ is a unit interval, $P_{3}$ is a hexagon, and $P_{4}$ is a truncated
octahedron. While $P_{4}$ was known to Johannes Kepler in 1619 \cite{JK}, the
family of permutahedra $P=\{P_{n}\}_{n\geq1}$ was first studied by the Dutch
mathematician Pieter Hendrik Schoute in 1911 \cite{PHS}.

The facets of $P_{n}$ are identified with the (ordered) \emph{partitions} of
$\underline{n}$ obtained by removing a single bar, called a \emph{divider},
from each permutation $A_{1}|\cdots|A_{n}$. For example, from the vertex
$3|1|2$ of $P_{3}$ we obtain the edges $13|2$ and $3|12$ with common vertex
$3|1|2$ (see Figure 1). In general, removing $n-k$ dividers from a permutation
of $\underline{n}$ produces a partition of length $k$ identified with an
$\left(  n-k\right)  $-dimensional face of $P_{n}$. \bigskip

\hspace{1.4in}\includegraphics[width=1.9in,height=1.50in]{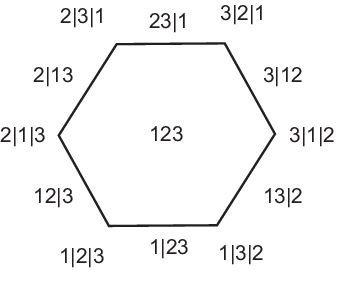} \medskip

\begin{center}
Figure 1. The Permutahedron $P_{3}$.\bigskip
\end{center}

There is a combinatorial bijection between the set of partitions of
$\underline{n}$ and the set of faces of $P_{n}$; the \emph{dimension} of a
partition of $\underline{n}$ is the geometric dimension of the corresponding
face. Furthermore, there is a combinatorial bijection between the set of faces
of $P_{n}$ and the set of \emph{down-rooted planar trees with levels and $n+1$
leaves} (DPTLs). Given a DPTL $T$ with $k$ levels and $n+1$ leaves, number the leaves
from left-to-right and assign the label $\ell$ to the vertex of $T$ at which
the branch containing leaf $\ell$ meets the branch containing leaf $\ell+1$ (a
vertex may have multiple labels). With level $1$ as the top level, let $A_{i}$
denote the set of vertex labels in level $i$; then $T$ corresponds to an
$(n-k)$-dimensional face $A_{1}|\cdots|A_{k}$ of $P_{n}$. For example, the DPTL
corresponding to the $2$-dimensional face $5|13|24$ of $P_{5}$ appears in
Figure 2. \vspace{0.5in}

\unitlength 1.5mm \linethickness{0.8pt}
\ifx\plotpoint\undefined\newsavebox{\plotpoint}\fi
\begin{picture}(77.572,16.187)(0,0)
	\put(36.5,18.443){\makebox(0,0)[cc]{$1\ \ \ \ 2\ \ \ \ 3\ \ \ 4\ \ \ 5\ \ \ 6$}}
	\put(43.8,12.4){\makebox(0,0)[cc]{$5$}}
	\put(31.0,11.0){\makebox(0,0)[cc]{$1$}}
	\put(35.2,11.0){\makebox(0,0)[cc]{$3$}}
	\put(36.2,7.0){\makebox(0,0)[cc]{$2\ \ 4$}}
	\multiput(28.019,15.977)(.070074384,-.067379215){117}{\line(1,0){.070074384}}
	\multiput(36.217,8.094)(.078612859,.067130082){119}{\line(1,0){.078612859}}
	\put(36.3,16.0){\line(0,-1){15.361}}
	\multiput(36.217,11.3)(.078612859,.067130082){43}{\line(1,0){.078612859}}
	\put(32.5,16.0){\line(0,-1){4.361}}
	\put(39.5,16.0){\line(0,-1){1.9}}
	\put(42.5,16.0){\line(0,-1){2.6}}
\end{picture}

\begin{center}
Figure 2. The DPTL corresponding to $5|13|24$.
\end{center}

\bigskip

Let $\theta^{n}$ denote the \emph{down-rooted} 
$n$\emph{-leaf corolla}, i.e., the DPTL with one level and $n$ leaves. 
One can extend a DPTL with $k$ levels to a DPTL with $k+1$ levels by grafting various corollas 
onto its leaves.  Thus the levels in a DPTL encode grafting order. 
For example, the DPTL in Figure 2 with one vertex in level 1, 
two vertices in level 2, and one vertex in level 3
is constructed by simultaneously grafting two copies of $\theta^{2}$ 
onto the first and second leaves of $\theta^{3}$  to form level 2, 
then grafting a third copy of $\theta^2$ onto the third leaf of $\theta^3$ to form level 1.

Now think of $\theta^{n}$ as a multilinear operation with one input and $n$ outputs, and 
denote the identity by $\mathbf{1}$.  A DPTL can be uniquely represented by a sequence of compositions
of $\theta^{n}$-operations and identities in the following way: Given a subsequence $i_1<\cdots<i_j$ of \underline{$m$} and operations $\theta^{n_1},\dots,\theta^{n_j}$, define the $\circ_{i_1\cdots i_j}$-composition 
\begin{equation}
\label{circ}
(\theta^{n_1},\dots,\theta^{n_j})\circ_{i_1 \dots i_j}\theta^{m}:=\left(  \mathbf{1}^{\otimes
i_1-1}\otimes\theta^{n_1}\cdots\otimes\mathbf{1}^{\otimes i_j-i_{j-1}-1}\otimes\theta^{n_j}\cdots\otimes \mathbf{1}^{m-i_j}\right)\theta^{m}  .
\end{equation}
Then the DPTL in Figure 2 is uniquely represented by the sequence of compositions
$$
\theta^2\circ_5((\theta^2,\theta^2)\circ_{1,2}\theta^3) = \left(\mathbf{1}^{\otimes4}\otimes\theta^{2}\right) 
\left(\theta^{2}\otimes\theta^{2}\otimes\mathbf{1}\right) \theta^{3}. 
$$
Equality of such compositions is an equivalence relation of the set of DPTLs. Thus $\left(  \theta_{2}\otimes\theta_{2}\otimes\theta_{2}\right)\theta_{3}  $
and  $ \left(\mathbf{1}^{\otimes4}\otimes\theta_{2}\right) 
\left(\theta_{2}\otimes\theta_{2}\otimes\mathbf{1}\right) \theta_{3} $ 
 are equivalent in $P_{5}.$

There is a related family of contractible polytopes $K=\{K_{n}\}_{n\geq2}$
called \emph{associahedra}, constructed by J. Stasheff in 1963 \cite{St}, whose faces
are indexed by up-rooted planar trees with $n$ leaves (without levels),
and there is a natural projection $\vartheta:P_{n-1}\rightarrow K_{n}$ due to
A. Tonks \cite{T} given by reflecting each DPTL with $n$ leaves in a horizontal axis and forgetting levels. 
Thus $\vartheta$ maps equivalent faces of $P_{n-1}$ to the same
face of $K_{n}.$

A \emph{non-}$\Sigma$\emph{ operad } $(\mathcal{P},\circ_{i})$ in the category
of sets consists of a sequence of sets $\mathcal{P}=\{\mathcal{P} 
(n)\}_{n\geq1}$, a unit map $\mathbf{1}\in\mathcal{P}(1)$, and products
$\circ_{i}:\mathcal{P}(m)\times\mathcal{P}(n)\rightarrow\mathcal{P}(m+n-1)$
for $m,n\geq1$ and $1\leq i\leq m$ interacting associatively \cite{MMS}.
Examples include (1) a collection of multilinear operations $\left\{
\theta^{n}:\theta^{1}=\mathbf{1}\right\}_{n\geq 2}  $ together with $\circ_{i_1} 
$-compositions, and (2) the top-dimensional cells of the associahedra $K$
identified with multilinear operations $\{\theta_{n}\}_{n\geq 2}$
(this time with $n$ inputs and one output) together with faces of $K$ identified with 
$\circ_{i_1}$-compositions of $\theta_n$-operations (the composition order in (\ref{circ}) is reversed). A \textit{map} 
of non-$\Sigma$ operads preserves operadic structure.

Stasheff defined the notion of an $A_{\infty}$\emph{-algebra} by appealing to
the operadic structure of the associahedra $K$ in the
following way: Denote the cellular chain complex of $K$ by $\left(  CC_{\ast
}\left(  K\right)  ,\partial\right)$.  Given a commutative ring $R$ with
unity and a differential graded (DG) $R$-module $\left(  A,d\right)  $, denote the tensor module of
$A$ by $TA$ and let $\nabla$ denote the differential on $Hom^{\ast}\left(
TA,A\right)  $ induced by $d,$ i.e., for $f\in Hom(A^{\otimes n},A)$ define
$$\nabla f: =  \sum_{i=0}^{n-1}f\left(\mathbf{1}^{\otimes i}\otimes
d\otimes\mathbf{1}^{\otimes n-i-1}\right) - (-1)^{|f|} df.$$

\begin{definition}
A DG $R$-module $(A,d)$ with differential $d$ of degree $+1$ together with a
family of multilinear operations $\{\omega_{n}\in Hom^{2-n}\left(  A^{\otimes
n},A\right) : n\geq 2 \}$ is an $A_{\infty}$\emph{-}\textbf{algebra} if there
exists a chain map of non-$\Sigma$ operads 
$\alpha:\left(  CC_{\ast}\left(  K\right)  ,\partial
\right)  \rightarrow \left(  Hom^{\ast}\left(  TA,A\right)  \},\nabla\right)  $
 such that $\alpha\left(  \theta_{n}\right)  =\omega_{n}$ for each $n.$
\end{definition}

S. Saneblidze and the second author extended 
the families of permutahedra and associahedra in \cite{SU1} and \cite{SU2}. 
The \emph{bipermutahedra} $PP=\{PP^n_m\}_{m+n\geq1}$ is a family of 
$\left(  m+n-1\right)  $-dimensional contractible polytopes, of which 
$PP^n_0\cong PP^0_n\cong P_n$.  The \emph{biassociahedra} 
$KK=\{KK^n_m\}_{mn\geq2}$ is a family of $\left(  m+n-3\right)  $-dimensional
contractible polytopes of which $KK^1_n\cong KK^n_1\cong K_n$. 
A picture of the heptagon $KK^2_3=PP_2^1$ appears in Figure 3. There is a natural
projection $\vartheta\vartheta:PP^n_m\rightarrow KK^{n+1}_{m+1}$, which agrees with
$\vartheta$ when $mn=0,$ and the notion of a \emph{matrad}, which generalizes 
non-$\Sigma$ operads (see \cite{SU1}). A \emph{map of matrads} preserves matradic structure. 
Identify the top-dimensional cell of $KK^n_m$ with a multilinear operation $\theta_{m}^{n}$ 
that has $m$ inputs and $n$ outputs.

\begin{definition}
A DG $R$-module $(A,d)$ with differential $d$ of degree $+1$ together with a
family of multilinear operations $\{\omega^n_m\in Hom^{3-m-n}\left(
A^{\otimes m},A^{\otimes n}\right):m+n\geq3 \}$ is an $A_{\infty}$\emph{-}
\hspace{-.05in}\textbf{bialgebra} if there exists a chain map of matrads $\alpha:\left(
CC_{\ast}\left(  KK\right)  ,\partial\right)  \rightarrow( Hom^{\ast}\left(
TA,TA\right)  ,\nabla) $ such that $\alpha\left(\theta_{m}^{n}\right)  =\omega^n_m$ 
for each $m$ and $n.$
\end{definition}

The definition of a matrad rests heavily on the notion of a \emph{bipartition
matrix }(BPM), whose entries are ordered pairs of partitions of the
same length (see Definitions \ref{defn2}) and whose rows
and columns are tightly controlled (see Definition \ref{defn3}). 
A BPM factors as a formal product of indecomposable BPMs, and the
dimension of a BPM is the sum of the dimensions of its indecomposable factors.
The dimension of an indecomposable BPM is the sum of its row, column, and 
entry dimensions (see Definition \ref{dim-GBPM}).

The faces of $PP^n_m$ are identified with certain products of
\textit{coherent generalized} BPMs (see \cite{SU2}). Just as the dimension of a partition of $\underline{n}$ is the
geometric dimension of the corresponding face of $P_{n}$, the dimension of a
product of coherent generalized BPMs corresponding to a face of $PP^n_m$ is the geometric
dimension of the face (see Figure 3). Thus computing the dimension of
generalized BPMs (and BPMs in particular) is fundamentally important.

The row, column, and entry dimensions of an indecomposable BPM are defined recursively 
and well-suited for machine computation. 
This article presents a computer program that computes the
dimension of a BPM by applying four routines of independent interest: (1) a
routine that factors a bipartition as a product of indecomposable BPMs,
(2) a routine that computes the inverse and recovers the bipartition, (3) a
routine that factors a BPM as a product of indecomposables, and (4)
a routine that calculates the ``transpose-rotation'' (the column dimension  of a BPM is the
row dimension of its transpose-rotation).

\section{Bipartitions and bipartition matrices}

\subsection{Partitions and bipartitions}

An \emph{ordered set} is either the empty set or a strictly increasing set of
positive integers.

\begin{definition}
\label{defn1}Let $A$ be a finite ordered set. A \textbf{partition of} $A$ 
\textbf{of length} $r$, denoted $A_{1}|A_{2}|\cdots|A_{r}$, is an ordered partition 
of ordered subsets of $A$. Each subset $A_{i}$ is called a \textbf{block}; two adjacent
blocks are separated by a \textbf{divider}. An empty block 
is denoted by $0$.
\end{definition}

Definition \ref{defn1} differs from the standard definition of a partition in
which all blocks are non-empty 
and unordered. This difference is important---the algorithm
for computing the dimension of a BPM essentially reduces to counting empty
blocks. Dividers will play a pivotal role in the definition of an ``embedding
partition,'' which encodes the precise way in which 
a subset $A$ of an ordered set $B$ is embedded in $B$ (see Definition \ref{defn4}).

\begin{example}
\label{ex1}The partition $p:=12|0|0$ has length $3$ and contains three
disjoint subsets of $A=\{1,2\}$. The first block contains the elements $1$ and
$2$, and as a stylistic choice we will not separate them by commas. The second
and third blocks are empty. In the context of the computer program in Appendix A, the
partition $p$ is represented as
\end{example}

\begin{lstlisting}
		[[1, 2], [], []].
\end{lstlisting}

\begin{definition}
\label{defn2}A \textbf{bipartition} of length $r$ is an ordered pair of 
partitions of length $r$ 
$$b:=\left(  A_{1}|A_{2}|\cdots|A_{r},\ B_{1}|B_{2}|\cdots
|B_{r}\right) $$ displayed as the fraction $$b=\frac{B_{1}|B_{2}|\cdots|B_{r}
}{A_{1}|A_{2}|\cdots|A_{r}^{\mathstrut}}.$$ The \textbf{input set of} $b$ is
defined and denoted by $IS(b):=A_{1}\cup A_{2}\cup\cdots\cup A_{r}$; the \textbf{output set
of} $b$ is defined and denoted by $OS(b):=B_{1}\cup B_{2}\cup\cdots\cup B_{r}$. An
\textbf{elementary bipartition} is a bipartition of length $1$; a \textbf{null
bipartition} has the form $\frac{0|\cdots|0}{0|\cdots|0^{\mathstrut}}$.
\end{definition}

\noindent Note that the partitions $A_{1}|A_{2}|\cdots|A_{r}$ and $B_{1}|B_{2} 
|\cdots|B_{r}$ in Definition \ref{defn2} may partition different sets; the
only requirement is that they have the same length.

\begin{example}
\label{ex2}The bipartition $b:=\frac{56|0|0}{1|2|4^{\mathstrut}}$ has length
$3$, $IS(b)=\{1,2,4\}$, and $OS(b)=\{5,6\}$. In the context of the computer program in
Appendix A, the bipartition $b$ is represented as
\end{example}

\begin{lstlisting}
		[[[1], [2], [4]], [[5, 6], [], []]].
\end{lstlisting}
\bigskip

In the development that follows, a bipartition $q/p$ is thought of as a
multivariable operator whose denominator $p$ is a sequence of inputs and whose
numerator $q$ is a sequence of outputs.

\subsection{Bipartition matrices}

\begin{definition}
\label{defn3}An $m\times n$ BPM $C=(c_{ij})$ has
the form
\[
C= 
\begin{pmatrix}
\frac{B_{1}^{11}|B_{2}^{11}|\cdots|B_{r_{11}\mathstrut}^{11}}{A_{1}^{11} 
|A_{2}^{11}|\cdots|A_{r_{11}}^{11\mathstrut}} & \cdots & \frac{B_{1} 
^{1n}|B_{2}^{1n}|\cdots|B_{r_{1n}\mathstrut}^{1n}}{A_{1}^{1n}|A_{2} 
^{1n}|\cdots|A_{r_{1n}}^{1n\mathstrut}}\\
\vdots & \ddots & \vdots\\
\frac{B_{1}^{m1}|B_{2}^{m1}|\cdots|B_{r_{m1}\mathstrut}^{m1}}{A_{1}^{m1} 
|A_{2}^{m1}|\cdots|A_{r_{m1}}^{m1\mathstrut}} & \cdots & \frac{B_{1} 
^{mn}|B_{2}^{mn}|\cdots|B_{r_{mn}\mathstrut}^{mn}}{A_{1}^{mn}|A_{2} 
^{mn}|\cdots|A_{r_{mn}}^{mn\mathstrut}} 
\end{pmatrix},
\]
and the following conditions hold for $1\leq i\leq m$ and $1\leq j\leq n$:

\begin{enumerate}
\item[1.] All bipartitions in a given column have equal input sets.

\item[2.] All bipartitions in a given row have equal output sets.

\item[3.] Either $IS(c_{ij}) = \varnothing$ or $\min IS(c_{ij}) > \max
IS(c_{kj}$ for all $k < i$.

\item[4.] Either $OS(c_{ij})=\varnothing$ or $\min OS(c_{ij})>\max OS(c_{ik})$
for all $k<j$.
\end{enumerate}

\noindent A \textbf{null }BPM has null bipartition entries. An \textbf{elementary} 
BPM has elementary bipartition entries.
\end{definition}

\noindent Items (3) and (4) in Definition \ref{defn3} indicate that the elements of the output sets in a given column
increase with the rows and the elements of the input sets in a given row
increase with the columns.

\begin{example}
\label{ex3}In the context of the computer program in Appendix A, the BPM
\[
C= 
\begin{pmatrix}
\frac{1|2}{3|0^{\mathstrut}} & \frac{12}{4^{\mathstrut}}\smallskip\\
\frac{5|6|7|8}{3|0|0|0^{\mathstrut}} & \frac{5|67|8}{0|0|4^{\mathstrut}} 
\end{pmatrix}
\]
is represented as
\end{example}

\begin{lstlisting}[breaklines]
	[[[[[3], []], [[1], [2]]], [[[4]], [[1, 2]]]], [[[[3], [], [], []], [[5], [6], [7], [8]]], [[[], [], [4]], [[5], [6, 7], [8]]]]].
\end{lstlisting}

\begin{example}
	The bipermutahedron $PP^1_2$ and the biassociahedron $KK^2_3$ are
	identical heptagons (see Figure 3). Its faces are identified with products of
	\textit{coherent} bipartition matrices whose
	dimension is the dimension of the corresponding face (for the definition of coherence see \cite{SU1}).
	
	\[
		{\includegraphics[
		height=2.55in,
		width=2.45in
		]{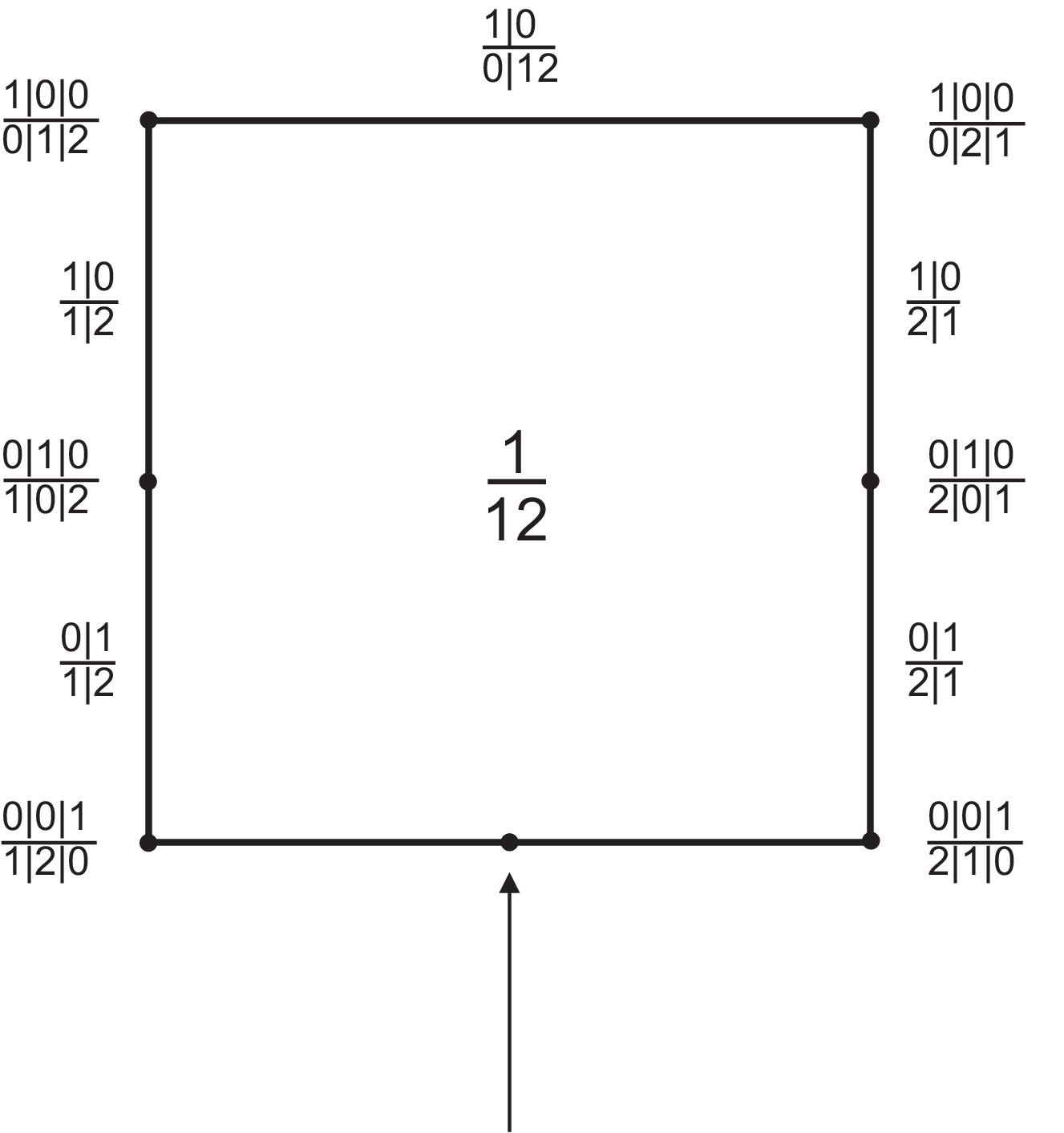}} 
	\]
	\vspace{-0.8in}
	\[
	\hspace{-0.15in}\left(
	\begin{array}
		[c]{c} 
		\frac{0}{12^{\mathstrut}}\medskip\\
		\frac{0|0}{1|2^{\mathstrut}} 
	\end{array}
	\right)  \left(
	\begin{array}
		[c]{ccc} 
		\frac{1}{0^{\mathstrut}} & \frac{1}{0^{\mathstrut}} & \frac{1}{0^{\mathstrut}} 
	\end{array}
	\right)  \text{\ \ \ \ \ }\left(
	\begin{array}
		[c]{c} 
		\frac{0|0}{2|1^{\mathstrut}}\medskip\\
		\frac{0}{12^{\mathstrut}} 
	\end{array}
	\right)  \left(
	\begin{array}
		[c]{ccc} 
		\frac{1}{0^{\mathstrut}} & \frac{1}{0^{\mathstrut}} & \frac{1}{0^{\mathstrut}} 
	\end{array}
	\right)
	\] 
	\[
	\left(
	\begin{array}
		[c]{c} 
		\frac{0|0}{2|1^{\mathstrut}}\medskip\\
		\frac{0|0}{1|2^{\mathstrut}} 
	\end{array}
	\right)  \left(
	\begin{array}
		[c]{ccc} 
		\frac{1}{0^{\mathstrut}} & \frac{1}{0^{\mathstrut}} & \frac{1}{0^{\mathstrut}} 
	\end{array}
	\right)
	\]
\end{example}
	\medskip
	\begin{center}
		Figure 3. The heptagon $KK^2_3=PP_2^1.$
	\end{center}

\medskip

\section{Embedding partitions}

Let $A=\{{a}_{1},{a}_{2},\dots,{a}_{m}\}$ and $B=\{{b}_{1}
,{b}_{2},\dots,{b}_{n}\}$ be ordered sets with $A\subseteq B$. The precise way in which $A$
embeds in $B$ is encoded by the \textquotedblleft embedding
partition.\textquotedblright  

\begin{definition}
\label{defn4}The \textbf{embedding partition of} $A$ \textbf{in} $B$, denoted by $\mathcal{EP} 
_{B}A$, is the partition of $A$ with the following properties:
\begin{enumerate}
\item[1.] The elements of $\mathcal{EP}_{B}A$ appear in the same order as in
$A$.

\item[2.] The number of dividers in $\mathcal{EP}_{B}A$
\begin{enumerate}
\item[$\bullet$] between consecutive elements ${a}_{i}$ and ${a}_{i+1}$ equals
the number of elements between ${a}_{i}$ and ${a}_{i+1}$ in $B$.

\item[$\bullet$] preceding ${a}_{1}$ equals the number of elements preceding
${a}_{1}$ in $B$.

\item[$\bullet$] succeeding ${a}_{m}$ equals the number of elements succeeding
${a}_{m}$ in $B$.

\item[$\bullet$] when $A=\varnothing$ equals the cardinality $\#B$.
\end{enumerate}
\end{enumerate}
\end{definition}

\noindent Note that the number of dividers in $\mathcal{EP}_{B}A$ is $\#B-\#A$.

\begin{remark}
	\noindent Live interactive Python implementations of the computational routines applied in the examples that follow are accessible at the link
	\begin{center}
		https://dawsonfreeman.pythonanywhere.com/
	\end{center}
\end{remark}

\begin{example}
\label{ex4} Let $A=\{3,5,8,9\}$ and $B=\{1,2,3,4,5,6,7,8,9\};$ then
\[
\mathcal{EP}_{B}A=0|0|3|5|0|89
\]
and has $\#B-\#A=5$ dividers. To compute $\mathcal{EP}_{B}A$ using the code in
Appendix A, input
\end{example}

\begin{lstlisting}[breaklines]
	print(embeddingPartition([1, 2, 3, 4, 5, 6, 7, 8, 9], [3, 5, 8, 9]));
\end{lstlisting}
\textit{the program returns} \begin{lstlisting}[breaklines]
	[[], [], [3], [5], [], [8, 9]].
\end{lstlisting}

\section{Factoring bipartitions}

When calculating BPM dimension, we first factor the BPM 
as a formal product of indecomposables.  Algorithm \ref{algor1} (below) factors
a $1 \times 1$ BPM $C$. Although $C$ does not equal its bipartition entry, we abuse notation and represent
$C$ as a bipartition and interpret the expression \textquotedblleft factoring a
bipartition\textquotedblright\ to mean \textquotedblleft factoring a
$1\times1$ BPM.\textquotedblright

\begin{algorithm}
\label{algor1} Let $c=\frac{B_{1}|B_{2}|\cdots|B_{r}}{A_{1}|A_{2}|\cdots
|A_{r}}$ be a bipartition with $r>1$.

\indent For $k=1$ to $r:$

\indent\indent Let $a_{k}^{1}|a_{k}^{2}|\cdots|a_{k}^{p_{k}}:=\mathcal{EP} 
_{A_{1}\cup\cdots\cup A_{k}}A_{k}$.

\indent\indent Let $b_{k}^{1}|b_{k}^{2}|\cdots|b_{k}^{q_{k}}:=\mathcal{EP} 
_{B_{k}\cup\cdots\cup B_{r}}B_{k}$.

\indent \indent Form the matrix
\[
C_{k} =
\begin{pmatrix}
\frac{b^{1}_{k_{\mathstrut}}}{a^{1^{\mathstrut}}_{k}} & \cdots & \frac
{b^{1}_{k_{\mathstrut}}}{a^{p^{\mathstrut}_{k}}_{k}}\\
\vdots & \ddots & \vdots\\
\frac{b^{q_{k}}_{k_{\mathstrut}}}{a^{1^{\mathstrut}}_{k_{\mathstrut}}} &
\cdots & \frac{b^{q_{k}}_{k_{\mathstrut}}}{a^{p^{\mathstrut}_{k}}_{k}} 
\end{pmatrix}
.
\]
\medskip

The \textbf{indecomposable factorization} of $c$ is the formal matrix product
\[
c=C_{1}\cdots C_{r}.
\]

\end{algorithm}

\noindent Clearly, the formal matrix product given by Algorithm \ref{algor1} is not
the standard matrix product in linear algebra.

\begin{example}
\label{ex5}The bipartition $\frac{0|34|5}{1|0|2^{\mathstrut}}$ factors as a
product of three BPMs whose entries are determined by
\[ 
\begin{matrix}
a_{1}^{1}=1 & b_{1}^{1}|b_{1}^{2}|b_{1}^{3}|b_{1}^{4}=0|0|0|0\\
a_{2}^{1}|a_{2}^{2}=0|0 & b_{2}^{1}|b_{2}^{2}=34|0\\
a_{3}^{1}|a_{3}^{2}=0|2 & b_{3}^{1}=5.
\end{matrix}
\]
Thus
\[
\frac{0|34|5}{1|0|2^{\mathstrut}}= 
\begin{pmatrix}
\frac{0}{1^{\mathstrut}}\medskip\\
\frac{0}{1^{\mathstrut}}\medskip\\
\frac{0}{1^{\mathstrut}}\medskip\\
\frac{0}{1^{\mathstrut}} 
\end{pmatrix} 
\begin{pmatrix}
\frac{34}{0^{\mathstrut}} & \frac{34}{0^{\mathstrut}}\medskip\\
\frac{0}{0^{\mathstrut}} & \frac{0}{0^{\mathstrut}} 
\end{pmatrix} 
\begin{pmatrix}
\frac{5}{0^{\mathstrut}} & \frac{5}{2^{\mathstrut}} 
\end{pmatrix}
.
\]
\medskip

\noindent To compute this factorization using the code in Appendix A, input
\end{example}

\begin{lstlisting}[breaklines]
	print(indecomposableFactorization([[[1], [], [2]], [[], [3, 4], [5]]])).
\end{lstlisting}
\medskip

\noindent\textit{The program returns the following list of bipartition
matrices:} \begin{lstlisting}[breaklines]
	[[[[[[1]], [[]]]], [[[[1]], [[]]]], [[[[1]], [[]]]], [[[[1]], [[]]]]], [[[[[]], [[3, 4]]], [[[]], [[3, 4]]]], [[[[]], [[]]], [[[]], [[]]]]], [[[[[]], [[5]]], [[[2]], [[5]]]]]].
\end{lstlisting}

\section{Equalizers}

Equalizers enable us to factor a BPM as a formal product of indecomposables and
play a pivotal role in the calculation of its dimension.  

\subsection{Row equalizers}

A \textquotedblleft row equalizer\textquotedblright\ transforms a BPM into a
BPM with equal numerators in each row.

\begin{definition}
\label{defn5} Let $C=(c_{ij})$ be a $q\times p$ BPM, where 

\[
c_{ij}:=\frac{B_{1}^{ij}|B_{2}^{ij}|\cdots|B_{r_{ij}\mathstrut}^{ij}} 
{A_{1}^{ij}|A_{2}^{ij}|\cdots|A_{r_{ij}}^{ij^{\mathstrut}}}\ .
\]
A \textbf{row equalizer of} $C$ is a $q\times p$ matrix of ordered sets
$E^{req}=(e_{ij}^{req}=\{e_{ij}^{1},\dots,e_{ij}^{n_{ij}}\})$ such that, when
the blocks of the corresponding bipartitions are combined in the following
way, the resulting output partitions in each row are equal:
\[
c_{ij}^{req}=\frac{(B_{1}^{ij}\cup\cdots\cup B_{e_{ij}^{1}}^{ij} 
)|({B_{e_{ij}^{1}+1}^{ij}}\cup\cdots\cup B_{e_{ij}^{2}}^{ij})|\cdots
|(B_{e_{ij}^{n}+1}^{ij}\cup\cdots\cup B_{r_{ij}}^{ij})}{(A_{1}^{ij}\cup
\cdots\cup A_{e_{ij}^{1}}^{ij})|(A_{e_{ij}^{1}+1}^{ij}\cup\cdots\cup
A_{e_{ij}^{2}}^{ij})|\cdots|(A_{e_{ij}^{n}+1}^{ij}\cup\cdots\cup A_{r_{ij} 
}^{ij})}\ .
\]
The matrix $C^{req}=\left(  c_{ij}^{req}\right)  $ is the corresponding
\textbf{row equalization of }$C.$\textbf{ }
\end{definition}

\noindent A row equalizer combines blocks in numerators and denominators so
that all numerators in the same row are equal. For example, a row equalizer
entry $e_{ij}^{req}=\{2,5\}$ indicates that blocks 3, 4, and 5 in the
numerator and denominator of $c_{ij}$ are to be combined.

\begin{example}
\label{ex6}Consider the matrix in Example \ref{ex3}:
\[
C= 
\begin{pmatrix}
\frac{1|2}{3|0^{\mathstrut}} & \frac{12}{4^{\mathstrut}}\medskip\\
\frac{5|6|7|8}{3|0|0|0^{\mathstrut}} & \frac{5|67|8}{0|0|4^{\mathstrut}} 
\end{pmatrix}
.
\]
To equalize numerators in the first row, combine all blocks in the numerator
and denominator of the first entry; then the first row becomes
\[
\frac{12}{3^{\mathstrut}}\ \ \ \frac{12}{4^{\mathstrut}}.
\]
\noindent To equalize numerators in the second row, combine second and third
blocks of the numerator and denominator of the first entry; then the second
row becomes
\[
\frac{5|67|8}{3|0|0^{\mathstrut}}\ \ \ \frac{5|67|8}{0|0|4^{\mathstrut}}.
\]
\noindent A row equalizer for $C$ is
\[
E^{req}= 
\begin{pmatrix}
\{\} & \{\}\smallskip\\
\{1,3\} & \{1,2\}
\end{pmatrix}
\]
and the corresponding row equalization is 
\[
C^{req}= 
\begin{pmatrix}
\frac{12}{3^{\mathstrut}} & \frac{12}{4^{\mathstrut}}\medskip\\
\frac{5|67|8}{3|0|0^{\mathstrut}} & \frac{5|67|8}{0|0|4^{\mathstrut}} 
\end{pmatrix}
.
\]

\end{example}

Note that an empty set in an equalizer indicates that all blocks in the
numerator and denominator of the corresponding entry are to be combined. As a
general rule, \emph{if a bipartition entry $c_{ij}$ has length $r_{ij}$ and
$e_{ij}^{req}\neq\varnothing$, then $\max e_{ij}^{req}<r_{ij}$.}

\subsection{Column equalizers}

A column equalizer transforms a BPM into a BPM with equal denominators in each column.

\begin{definition}
\label{defn6}Let $C=(c_{ij})$ be a $q\times p$ BPM, where 

\[
c_{ij}=\frac{B_{1}^{ij}|B_{2}^{ij}|\cdots|B_{r_{ij}}^{ij}}{A_{1}^{ij} 
|A_{2}^{ij}|\cdots|A_{r_{ij}}^{ij^{\mathstrut}}}\ .
\]
A \textbf{column equalizer of} $C$ is a $q\times p$ matrix of ordered sets
$E^{ceq}=(e_{ij}^{ceq}=\{e_{ij}^{1},\dots,e_{ij}^{n_{ij}}\})$ such that, when
the blocks of the corresponding bipartition are combined in the following way,
the resulting input partitions in each column are equal:
\[
c_{ij}^{ceq}=\frac{(B_{1}^{ij}\cup\cdots\cup B_{e_{ij}^{1}}^{ij} 
)|({B_{e_{ij}^{1}+1}^{ij}}\cup\cdots\cup B_{e_{ij}^{2}}^{ij})|\cdots
|(B_{e_{ij}^{n}+1}^{ij}\cup\cdots\cup B_{r_{ij}}^{ij})}{(A_{1}^{ij}\cup
\cdots\cup A_{e_{ij}^{1}}^{ij})|(A_{e_{ij}^{1}+1}^{ij}\cup\cdots\cup
A_{e_{ij}^{2}}^{ij})|\cdots|(A_{e_{ij}^{n}+1}^{ij}\cup\cdots\cup A_{r_{ij} 
}^{ij})}\ .
\]
The matrix $C^{ceq}=\left(  c_{ij}^{ceq}\right)  $ is the \textbf{column
equalization of }$C.$\textbf{ }
\end{definition}

\begin{example}
\label{ex7} Consider the matrix in Example \ref{ex3}:
\[
C= 
\begin{pmatrix}
\frac{1|2}{3|0^{\mathstrut}} & \frac{12}{4^{\mathstrut}}\medskip\\
\frac{5|6|7|8}{3|0|0|0^{\mathstrut}} & \frac{5|67|8}{0|0|4^{\mathstrut}} 
\end{pmatrix}
.
\]
To equalize denominators in the first column, combine second, third, and
fourth blocks in the numerator and denominator of the entry in the second row;
then the first column becomes
\[
\left(  \frac{1|2}{3|0^{\mathstrut}}\ \ \ \frac{5|678}{3|0^{\mathstrut} 
}\right)  ^{T}.
\]
To equalize denominators in the second column, combine all blocks in the
numerator and denominator of the entry in the second row; then the second
column becomes
\[
\left(  \frac{12}{4^{\mathstrut}}\ \ \ \frac{5678}{4^{\mathstrut}}\right)
^{T}.
\]
A column equalizer for the matrix $C$ is
\[
E^{ceq}= 
\begin{pmatrix}
\{1\} & \{\}\smallskip\\
\{1\} & \{\}
\end{pmatrix}
\]
and the corresponding equalization is
\[
C^{ceq}= 
\begin{pmatrix}
\frac{1|2}{3|0^{\mathstrut}} & \frac{12}{4^{\mathstrut}}\medskip\\
\frac{5|678}{3|0^{\mathstrut}} & \frac{5678}{4^{\mathstrut}} 
\end{pmatrix}
.
\]

\end{example}

\subsection{Equalizers}

\begin{definition}
\label{defn7}An \textbf{equalizer} of a $q\times p$ BPM $C$ is a $q\times p$
matrix of ordered sets $E$ that is simultaneously a row equalizer and a column
equalizer of $C$.
\end{definition}

\noindent The program in Appendix A generates a list of all row equalizers of
$C$ and a list of all column equalizers of $C$. Matrices in their intersection are
the equalizers of $C$.

\begin{example}
\label{ex8}Consider the BPM $C$ in Example \ref{ex3}: 
\[
C= 
\begin{pmatrix}
\frac{1|2}{3|0^{\mathstrut}} & \frac{12}{4^{\mathstrut}}\medskip\\
\frac{5|6|7|8}{3|0|0|0^{\mathstrut}} & \frac{5|67|8}{0|0|4^{\mathstrut}} 
\end{pmatrix}
.
\]
The only way to equalize numerators in each row and denominators in each
column is to combine all blocks. The equalized matrix is the elementary
matrix
\[
C^{eq}= 
\begin{pmatrix}
\frac{12}{3^{\mathstrut}} & \frac{12}{4^{\mathstrut}}\medskip\\
\frac{5678}{3^{\mathstrut}} & \frac{5678}{4^{\mathstrut}} 
\end{pmatrix}
.
\]
Its equalizer is the \textbf{null equalizer}
\[
E_{\varnothing}= 
\begin{pmatrix}
\{\} & \{\}\smallskip\\
\{\} & \{\}
\end{pmatrix}
.
\]

\end{example}

\noindent One can always trivially equalize a BPM by combining all blocks as
in Example \ref{ex8}.

\begin{example}
\label{ex9}A non-trivial equalization of the BPM
\[
C= 
\begin{pmatrix}
\frac{1|2|0}{3|0|4} & \frac{1|2}{5|0}\smallskip\\
\frac{6|7|8|9}{0|3|0|4} & \frac{6|7|8|9}{0|5|0|0} 
\end{pmatrix}
\text{\ \ is\ \ }
C^{eq}= 
\begin{pmatrix}
\frac{1|2}{3|4} & \frac{1|2}{5|0}\smallskip\\
\frac{67|89}{3|4} & \frac{67|89}{5|0} 
\end{pmatrix}
\]
with equalizer
\[
E= 
\begin{pmatrix}
\{1\} & \{1\}\smallskip\\
\{2\} & \{2\}
\end{pmatrix}
.
\]

\end{example}

\begin{definition}
\label{defn8}Let $C$ be a $q\times p$ BPM. Let $\mathcal{EP}$ denote the set
of all equalizers of $C$ with the following property: The cardinality of each
entry is greater than or equal to the cardinalities of the corresponding
entries in all equalizers of $C$. The \textbf{maximal equalizer} of $C$ is the
matrix in $\mathcal{EP}$ whose entries are less than the corresponding entries
in all other matrices in $\mathcal{EP}$ with respect to lexicographic order.
\end{definition}

\begin{example}
To compute the maximal equalizer of the matrix in Example \ref{ex3} using the
code in Appendix A, input
\end{example}

\begin{lstlisting}[breaklines]
	print(getMaximalEqualizer([[[[[3], []], [[1], [2]]], [[[4]], [[1, 2]]]], [[[[3], [], [], []], [[5], [6], [7], [8]]], [[[], [], [4]], [[5], [6, 7], [8]]]]])).
	
\end{lstlisting}
\textit{The program returns the maximal equalizer}
\begin{lstlisting}[breaklines]
	[[[], []], [[], []]].
\end{lstlisting}

\begin{example}
\label{ex10}Suppose a $2\times2$ BPM\ $C$ has the following equalizers:
\[
E_{1}=
\begin{pmatrix}
\{1,2\} & \{2,3\}\smallskip\\
\{1,3\} & \{1,2\}
\end{pmatrix}
;\ \ E_{2}=
\begin{pmatrix}
\{1\} & \{2\}\smallskip\\
\{2\} & \{1\}
\end{pmatrix}
;
\]
\[
E_{3}=
\begin{pmatrix}
\{1,2\} & \{1,3\}\smallskip\\
\{1,3\} & \{1,2\}
\end{pmatrix}
;\ \ E_{4}=
\begin{pmatrix}
\{\} & \{\}\smallskip\\
\{\} & \{\}
\end{pmatrix}
.
\]
The entries of $E_{1}$ and $E_{3}$ have cardinality $2$, the entries of
$E_{2}$ have cardinality $1$, and the entries of $E_{4}$ have cardinality $0$.
Since $E_{3}$ precedes $E_{1}$ in a least-to-greatest lexicographically sorted
list, $E_{3}$ is the maximal equalizer of $C$.
\end{example}

\begin{example}
\label{ex11}Refer to Example \ref{ex9} and consider the equalizers
\[
E_{\varnothing}= 
\begin{pmatrix}
\{\} & \{\}\smallskip\\
\{\} & \{\}
\end{pmatrix}
\text{ \ and \ }E= 
\begin{pmatrix}
\{1\} & \{1\}\smallskip\\
\{2\} & \{2\}
\end{pmatrix}
.
\]
The maximal equalizer is
\[ 
\begin{pmatrix}
\{1\} & \{1\}\smallskip\\
\{2\} & \{2\}
\end{pmatrix}
.
\]

\end{example}

\section{The transverse decomposition}

As mentioned above, the first step in computing BPM dimension is to factor the BPM 
as a product of indecomposables. The \textquotedblleft transverse 
decomposition\textquotedblright\ of a decomposable BPM $C$ is a
factorization of the form $C=AB$, where $A$ is an indecomposable column matrix and $B$
is a (possibly decomposable) row matrix.

\begin{algorithm}
\label{algor2}Let $C=\left(  c_{ij}\right)  $ be an $m\times n$ BPM with
non-null maximal equalizer $E=(e_{ij})$.

\indent For each $c_{ij}=\frac{B_{1}^{ij}|B_{2}^{ij}|\cdots|B_{r_{ij}}^{ij} 
}{A_{1}^{ij}|A_{2}^{ij}|\cdots|A_{r_{ij}}^{ij^{\mathstrut}}}$ in $C:$

\indent\indent Let $e$ be the first element of $e_{ij}$.

\indent\indent Let $a_{1}^{ij}|\cdots|a_{p}^{ij}:=\mathcal{EP}_{IS(c_{ij} 
)}(A_{e+1}^{ij}\cup\cdots\cup A_{r_{ij}}^{ij})$.

\indent\indent For $k$ from $1$ to $p:$

\indent\indent\indent Let $a_{k,e+1}^{ij}|\cdots|a_{k,r_{ij}}^{ij} 
:=(A_{e+1}^{ij}\cap a_{k}^{ij})|\cdots|(A_{r_{ij}}^{ij}\cap a_{k}^{ij})$.
\medskip

\indent \indent Let $b^{ij}_{1} | \cdots| b^{ij}_{q} := \mathcal{EP}
_{OS(c_{ij})}(B^{ij}_{1} \cup\cdots\cup B^{ij}_{e})$.

\indent\indent For $k$ from $1$ to $q:$

\indent\indent\indent Let $b_{k,1}^{ij}|\cdots|b_{k,e}^{ij}:=(B_{1}^{ij}\cap
b_{k}^{ij})|\cdots|(B_{e}^{ij}\cap b_{k}^{ij})$.
\medskip

\indent\indent Define the formal matrix decomposition
\[
C_{1}^{ij}C_{2}^{ij}= 
\begin{pmatrix}
\frac{b_{1,1}^{ij}|\cdots|b_{1,e}^{ij}}{A_{1}^{ij}|\cdots|A_{e} 
^{ij^{\mathstrut}}}\\
\vdots\\
\frac{b_{q,1}^{ij}|\cdots|b_{q,e}^{ij}}{A_{1}^{ij}|\cdots|A_{e} 
^{ij^{\mathstrut}}} 
\end{pmatrix} 
\begin{pmatrix}
\frac{B_{e+1}^{ij}|\cdots|B_{r_{ij}}^{ij}}{a_{1,e+1}^{ij}|\cdots|a_{1,r_{ij} 
}^{ij^{\mathstrut}}} & \cdots & \frac{B_{e+1}^{ij}|\cdots|B_{r_{ij}}^{ij} 
}{a_{p,e+1}^{ij}|\cdots|a_{p,r_{ij}}^{ij^{\mathstrut}}} 
\end{pmatrix}
.
\]
The \textbf{transverse decomposition of} $C$ \textbf{with respect to} $E$ is
the formal matrix factorization $C:=C_{1}C_{2}$ where
\[
C_{1}=\left(  C_{1}^{ij}\right)  = 
\begin{pmatrix}
\frac{b_{1,1}^{11}|\cdots|b_{1,e}^{11}}{A_{1}^{11}|\cdots|A_{e}
^{11^{\mathstrut}}} &  & \frac{b_{1,1}^{1n}|\cdots|b_{1,e}^{1n}}{A_{1}
^{1n}|\cdots|A_{e}^{1n^{\mathstrut}}}\\
\vdots & \cdots & \vdots\\
\frac{b_{q,1}^{11}|\cdots|b_{q,e}^{11}}{A_{1}^{11}|\cdots|A_{e}
^{11^{\mathstrut}}} &  & \frac{b_{q,1}^{1n}|\cdots|b_{q,e}^{1n}}{A_{1}
^{1n}|\cdots|A_{e}^{1n^{\mathstrut}}}\\
\vdots &  & \vdots\\
\frac{b_{1,1}^{m1}|\cdots|b_{1,e}^{m1}}{A_{1}^{m1}|\cdots|A_{e}
^{m1^{\mathstrut}}} &  & \frac{b_{1,1}^{mn}|\cdots|b_{1,e}^{mn}}{A_{1}
^{mn}|\cdots|A_{e}^{mn^{\mathstrut}}}\\
\vdots & \cdots & \vdots\\
\frac{b_{q,1}^{m1}|\cdots|b_{q,e}^{m1}}{A_{1}^{m1}|\cdots|A_{e}
^{m1^{\mathstrut}}} &  & \frac{b_{q,1}^{mn}|\cdots|b_{q,e}^{mn}}{A_{1}
^{mn}|\cdots|A_{e}^{mn^{\mathstrut}}} 
\end{pmatrix}
\]
and $C_{2}=\left(  C_{2}^{ij}\right)  =$
\[ 
\begin{pmatrix}
\frac{B_{e+1}^{11}|\cdots|B_{r_{ij}}^{11}}{a_{1,e+1}^{11}|\cdots|a_{1,r_{ij}
}^{11^{\mathstrut}}} & \cdots & \frac{B_{e+1}^{11}|\cdots|B_{r_{ij}}^{11}
}{a_{p,e+1}^{11}|\cdots|a_{p,r_{ij}}^{11^{\mathstrut}}} & \cdots &
\frac{B_{e+1}^{1n}|\cdots|B_{r_{ij}}^{1n}}{a_{1,e+1}^{1n}|\cdots|a_{1,r_{ij}
}^{1n^{\mathstrut}}} & \cdots & \frac{B_{e+1}^{1n}|\cdots|B_{r_{ij}}^{1n}
}{a_{p,e+1}^{1n}|\cdots|a_{p,r_{ij}}^{1n^{\mathstrut}}}\\
& \vdots &  &  &  & \vdots & \\
\frac{B_{e+1}^{m1}|\cdots|B_{r}^{m1}}{a_{1,e+1}^{m1}|\cdots|a_{1,r}
^{m1^{\mathstrut}}} & \cdots & \frac{B_{e+1}^{m1}|\cdots|B_{r_{ij}}^{m1}
}{a_{p,e+1}^{m1}|\cdots|a_{p,r_{ij}}^{m1^{\mathstrut}}} & \cdots &
\frac{B_{e+1}^{mn}|\cdots|B_{r}^{mn}}{a_{1,e+1}^{mn}|\cdots|a_{1,r_{ij}
}^{mn^{\mathstrut}}} & \cdots & \frac{B_{e+1}^{mn}|\cdots|B_{r}^{mn}
}{a_{p,e+1}^{mn}|\cdots|a_{p,r_{ij}}^{mn^{\mathstrut}}} 
\end{pmatrix}
.
\]
\end{algorithm}

\noindent The entries $C_{1}^{ij}$ and $C_{2}^{ij}$ have the following properties:

\begin{enumerate}
\item[1.] $C_{1}^{ij}$ is indecomposable.

\item[2.] $1+\#IS(C_{1}^{ij})$ equals the number of entries in $C_{2}^{ij}$.

\item[3.] $1+\#OS(C_{2}^{ij})$ equals the number of entries in $C_{1}^{ij}$.
\end{enumerate}
\medskip

\begin{example}
To factor the (indecomposable) matrix in Example \ref{ex8} using the code in
Appendix A, input
\end{example}

\begin{lstlisting}[breaklines]
	print(factorization([[[[[3], []], [[1], [2]]], [[[4]], [[1, 2]]]], [[[[3], [], [], []], [[5], [6], [7], [8]]], [[[], [], [4]], [[5], [6, 7], [8]]]]]))
\end{lstlisting}
\textit{Since the matrix is indecomposable, the program
returns the input matrix.} \medskip

\begin{example}
The (unique) indecomposable factorization of
\[
C=
\begin{pmatrix}
\frac{2|3|0|4}{1|3|2|0^{\mathstrut}} &  & \frac{23|0|4}{7|5|6^{\mathstrut}}\\
\frac{7|56|0}{13|0|2^{\mathstrut}} &  & \frac{7|0|56}{5|7|6^{\mathstrut}} 
\end{pmatrix}
\]
is 
\[
C_{1}C_{2}=\left(
\begin{array}
[c]{cc} 
\frac{2|3}{1|3^{\mathstrut}} & \frac{23|0}{7|5^{\mathstrut}}\\
\frac{0|0}{1|3^{\mathstrut}} & \frac{0|0}{7|5^{\mathstrut}}\\
\frac{0}{13^{\mathstrut}} & \frac{0|0}{5|7^{\mathstrut}}\\
\frac{0}{13^{\mathstrut}} & \frac{0|0}{5|7^{\mathstrut}}\\
\frac{7}{13^{\mathstrut}} & \frac{7|0}{5|7^{\mathstrut}} 
\end{array}
\right)  \left(
\begin{array}
[c]{cccccc} 
\frac{0|4}{0|0^{\mathstrut}} & \frac{0|4}{2|0^{\mathstrut}} & \frac
{0|4}{0|0^{\mathstrut}} & \frac{4}{0^{\mathstrut}} & \frac{4}{6^{\mathstrut}}
& \frac{4}{0^{\mathstrut}}\\
\frac{56|0}{0|0^{\mathstrut}} & \frac{56|0}{0|2^{\mathstrut}} & \frac
{56|0}{0|0^{\mathstrut}} & \frac{56}{0^{\mathstrut}} & \frac{56} 
{6^{\mathstrut}} & \frac{56}{0^{\mathstrut}} 
\end{array}
\right)  .
\]
\end{example}
\medskip
\noindent \textit{Using the code in Appendix A, input}

\begin{lstlisting}[breaklines]
print(factorization([[[[[1],[3],[2],[]],[[2],[3],[],[4]]],[[[7],[5],[6]],
      [[2,3],[],[4]]]],[[[[1,3],[],[2]],[[7],[5,6],[]]],
      [[[5],[7],[6]],[[7],[],[5,6]]]]]))
\end{lstlisting}
\textit{The program returns the indecomposable factorization}
\begin{lstlisting}[breaklines]
[[[[[[1], [3]], [[2], [3]]], [[[7], [5]], [[2, 3], []]]], [[[[1], [3]], [[], []]], [[[7], [5]], [[], []]]], [[[[1, 3]], [[]]], [[[5], [7]], [[], []]]], [[[[1, 3]], [[]]], [[[5], [7]], [[], []]]], [[[[1, 3]], [[7]]], [[[5], [7]], [[7], []]]]], [[[[[], []], [[], [4]]], [[[2], []], [[], [4]]], [[[], []], [[], [4]]], [[[]], [[4]]], [[[6]], [[4]]], [[[]], [[4]]]], [[[[], []], [[5, 6], []]], [[[], [2]], [[5, 6], []]], [[[], []], [[5, 6], []]], [[[]], [[5, 6]]], [[[6]], [[5, 6]]], [[[]], [[5, 6]]]]]] .
\end{lstlisting}

Given an $m\times n$ BPM  $C=\left(  c_{ij}\right)  $, compute the indecomposable 
factorization of $C$ in the following way: If the maximal
equalizer $E_{1}$ of $C$ is null, $C=C_{1}$ is the indecomposable
factorization. Otherwise, compute the transverse decomposition $C=C_{1}C_{2}$
with respect to $E_{1}$. If the maximal equalizer $E_{2}$ of $C_{2}$ is null,
$C=C_{1}C_{2}$ is the indecomposable factorization. Otherwise, compute the
transverse decomposition $C_{2}:=C_{2}^{\prime}C_{3}$ with respect to $E_{2}$.
If the maximal equalizer of $C_{3}$ is null, the indecomposable factorization
is $C=C_{1}C_{2}^{\prime}C_{3}$. Otherwise, continue this process for $s-1$
steps until the maximal equalizer $E_{s}$ of $C_{s}$ is null. Then
$C=C_{1}C_{2}^{\prime}\cdots C_{s-1}^{\prime}C_{s}$ is the indecomposable factorization.

The tools we need to compute the dimension of a BPM are now in place.

\section{The dimension algorithm}

Given a BPM $C$, let $\pi\left(  C\right)  $ denote the matrix obtained from
$C$ by discarding empty biblocks $0/0$ in non-null entries and by reducing
null entries to $0/0$. Denote the dimension of $C$ by
$\left\vert C\right\vert .$ 

\begin{definition}
\label{dim-GBPM}Let $C$ be a BPM. If $C$ is null, define $\left\vert C\right\vert
:=0.$ Otherwise, define
\begin{equation}
\left\vert C\right\vert :=\left\vert C\right\vert ^{row}+\left\vert
C\right\vert ^{col}+|C|^{ent}, \label{dim-framed} 
\end{equation}
where the \textbf{row dimension }$\left\vert C\right\vert ^{row}$\textbf{,
column dimension }$\left\vert C\right\vert ^{col}$, and \textbf{entry
dimension} $\left\vert C\right\vert ^{ent}$ are independent and given by the
following recursive algorithms:
\bigskip

\noindent\textbf{Row dimension algorithm. }

Given a $q\times p$ BPM $C$, let $C_{1}\cdots C_{r}$ be its indecomposable factorization.  \smallskip

\smallskip If $r>1,$ define $\left\vert C\right\vert ^{row} 
:=\sum_{k\in\underline{r}}\left\vert C_{k}\right\vert ^{row},$ 
\newline\hspace*{.4in}where
$\left\vert C_{k}\right\vert ^{row}$ is given by setting $C=C_{k}$ and continuing 
recursively.\smallskip

If $r=1$ and $q>1$, define $\left\vert C\right\vert ^{row} 
:=\sum_{i\in\underline{q}}\left\vert C_{i\ast}\right\vert ^{row},$ 
\newline\hspace*{.4in}where
$\left\vert C_{i\ast}\right\vert ^{row}$ is given by setting $C=C_{i\ast}$ and 
continuing recursively.\smallskip

Otherwise, $C$ is a row matrix $\left(  c_{1}\cdots c_{p}\right)  .$
\newline\hspace*{.4in}Set $C=\pi(C)$. \smallskip

\hspace*{.2in}If $C$ is an elementary matrix and $OS(C)\neq\varnothing,$ define $\left\vert
C\right\vert ^{row}:=0.$\medskip

\hspace*{.2in}If $C=\left(  \frac{0}{\mathbf{a}_{1}}\cdots\frac{0}{\mathbf{a}_{p}}\right)
,$ define $\left\vert C\right\vert ^{row}:=\left\{
\begin{array}
[c]{cc} 
0, & C\text{ is null}\\
\#IS(C)-1, & \text{otherwise.} 
\end{array}
\right.  $\medskip

\hspace*{.2in}Otherwise, define $\left\vert C\right\vert ^{row}:=\sum_{j\in\mathfrak{p} 
}\left\vert c_{j}\right\vert ^{row},$ 
\newline\hspace*{.6in}where $\left\vert c_{j}\right\vert
^{row}$ is given by setting $C=c_{j}$ and 
continuing recursively.\bigskip

\noindent\textbf{Column dimension algorithm. }

Given a $q\times p$ BPM $C$, let $C_{1}\cdots C_{r}$ be its indecomposable factorization. \smallskip

If $r>1,$ define $\left\vert C\right\vert ^{col}:=\sum
_{k\in\underline{r}}\left\vert C_{k}\right\vert ^{col},$ 
\newline\hspace*{.4in}where $\left\vert
C_{k}\right\vert ^{col}$ is given by setting $C=C_{k}$ and continuing 
recursively.\smallskip

If $r=1$ and $p>1$, define $\left\vert C\right\vert ^{col} 
:=\sum_{j\in\underline{p}}\left\vert C_{\ast j}\right\vert ^{col},$ 
\newline \hspace*{.4in}where
$\left\vert C_{\ast j}\right\vert ^{col}$ is given by setting  $C=C_{\ast j}$ and 
continuing recursively.\smallskip

Otherwise, $C$ is a column matrix $\left(  c_{1}\cdots c_{q}\right)
^{T}.$ 
\newline \hspace*{.4in}Set $C=\pi(C)$.\smallskip

\hspace*{.2in}If $C$ is an elementary matrix and $IS(C)\neq\varnothing$, define $\left\vert
C\right\vert ^{col}:=0.$\medskip

\hspace*{.2in}If $C=\left(  \frac{\mathbf{b}_{1}}{0}\cdots\frac{\mathbf{b}_{q}}{0}\right)
^{T},$ define $\left\vert C\right\vert ^{col}:=\left\{
\begin{array}
[c]{cc} 
0, & C\text{ is null}\\
\#OS(C)-1, & \text{otherwise.} 
\end{array}
\right.  $\medskip

\hspace*{.2in}Otherwise, define $|C|^{col}:=\sum_{i\in\underline{q}}\left\vert
c_{i}\right\vert ^{col},$ 
\newline\hspace*{.6in}where $\left\vert c_{i}\right\vert ^{col}$ is given
by setting $C=c_{i}$ and continuing 
recursively.\bigskip

\noindent\textbf{Entry dimension algorithm. }

Given a $q\times p$ BPM $C$, let $C_{1}\cdots C_{r}$ be its indecomposable factorization. \smallskip

If $r>1,$ define $\left\vert C\right\vert ^{ent}:=\sum
_{k\in\underline{r}}\left\vert C_{k}\right\vert ^{ent},$ 
\newline\hspace*{.4in}where $\left\vert
C_{k}\right\vert ^{ent}$ is given by setting $C=C_{k}$ and continuing recursively.\smallskip

Otherwise, define $|C|^{ent}:=\sum_{\left(  i,j\right)
\in\underline{q}\times\underline{p}}\left\vert c_{ij}\right\vert ^{ent},$
\newline\hspace*{.4in}where $\left\vert c_{ij}\right\vert ^{ent}$ is given by setting $C=c_{ij}$ and continuing recursively
\smallskip \newline\hspace*{.4in}unless $c_{ij}=\frac{\mathbf{b}_{i}}{\mathbf{a}_{j}},$ in which case define
\[
\left\vert c_{ij}\right\vert ^{ent}:=\left\{
\begin{array}
[c]{cl} 
\#\mathbf{a}_{j}+\#\mathbf{b}_{i}-1, & \mathbf{a}_{j},\mathbf{b}_{i} 
\neq\varnothing\\
0, & \text{otherwise.} 
\end{array}
\right.
\]
\end{definition}
\smallskip

\begin{example}
Consider the matrix in Example \ref{ex3}:
\[
C= 
\begin{pmatrix}
\frac{1|2}{3|0^{\mathstrut}} & \frac{12}{4^{\mathstrut}}\medskip\\
\frac{5|6|7|8}{3|0|0|0^{\mathstrut}} & \frac{5|67|8}{0|0|4^{\mathstrut}} 
\end{pmatrix}
.
\]
Using the program in Appendix A, the row dimension is $0$, the column dimension
is $1$, and the entry dimension is $5$. Therefore the dimension of $C$ is $6$.
A similar example with output given by the program in Appendix A appears in
Appendix B.
\end{example}

\begin{example}
\label{ex22} Consider the matrix
\[ 
\begin{pmatrix}
\frac{1|2|0}{0|1|2^{\mathstrut}} & \frac{1|0|2}{0|3|4^{\mathstrut}} &
\frac{1|2|0|0}{5|0|6|0^{\mathstrut}} & \frac{1|2}{7|8^{\mathstrut}}\medskip\\
\frac{3|4}{0|12^{\mathstrut}} & \frac{3|0|4}{0|3|4^{\mathstrut}} &
\frac{3|0|4|0}{5|0|6|0^{\mathstrut}} & \frac{3|4}{7|8^{\mathstrut}}\medskip\\
\frac{5|6|0}{0|1|2^{\mathstrut}} & \frac{5|6|0}{0|3|4^{\mathstrut}} &
\frac{5|6|0|0}{5|0|6|0^{\mathstrut}} & \frac{5|6}{7|8^{\mathstrut}}\medskip\\
\frac{7|8}{0|12^{\mathstrut}} & \frac{7|8}{0|34^{\mathstrut}} & \frac
{7|8|0|0}{5|0|6|0^{\mathstrut}} & \frac{7|8|0|0}{7|0|0|8^{\mathstrut}} 
\end{pmatrix}
.
\]
Using the program in Appendix A, the row dimension is $4$, the column dimension
is $10$, and the entry dimension is $23$. Therefore the dimension of $C$ is
$37$.
\end{example}

\begin{example}
Consider the matrix
\[ 
\begin{pmatrix}
\frac{1|2|0}{0|1|2^{\mathstrut}} & \frac{1|0|2}{0|3|4^{\mathstrut}} &
\frac{1|2|0|0}{5|0|6|0^{\mathstrut}} & \frac{1|2}{7|8^{\mathstrut}} &
\frac{1|2}{0|9^{\mathstrut}}\medskip\\
\frac{3|4}{0|12^{\mathstrut}} & \frac{3|0|4}{0|3|4^{\mathstrut}} &
\frac{3|0|4|0}{5|0|6|0^{\mathstrut}} & \frac{3|4}{7|8^{\mathstrut}} &
\frac{3|4|0}{0|9|0^{\mathstrut}}\medskip\\
\frac{5|6|0}{0|1|2^{\mathstrut}} & \frac{5|6|0}{0|3|4^{\mathstrut}} &
\frac{5|6|0|0}{5|0|6|0^{\mathstrut}} & \frac{5|6}{7|8^{\mathstrut}} &
\frac{5|6}{0|9^{\mathstrut}}\medskip\\
\frac{7|8}{0|12^{\mathstrut}} & \frac{7|8}{0|34^{\mathstrut}} & \frac
{7|8|0|0}{5|0|6|0^{\mathstrut}} & \frac{7|8|0|0}{7|0|0|8^{\mathstrut}} &
\frac{7|0|8}{0|9|0^{\mathstrut}}\medskip\\
\frac{9|0}{0|12^{\mathstrut}} & \frac{9|0}{0|34^{\mathstrut}} & \frac
{9|0}{5|6^{\mathstrut}} & \frac{9|0|0|0}{7|0|8|0^{\mathstrut}} & \frac
{9|0|0}{0|9|0^{\mathstrut}} 
\end{pmatrix}
.
\]
This matrix has one more row and one more column than the matrix in Example
\ref{ex22}. The row dimension of this matrix is $6$, the column dimension is
$12$, and the entry dimension is $28$. Therefore, its dimension is $46$.
\end{example}

\begin{example}
Consider the matrix
\[ 
\begin{pmatrix}
\frac{1|2|0}{0|1|2^{\mathstrut}} & \frac{1|0|2}{0|3|4^{\mathstrut}} &
\frac{1|2|0|0}{5|0|6|0^{\mathstrut}} & \frac{1|2}{7|8^{\mathstrut}}\medskip\\
\frac{3|4}{12|0^{\mathstrut}} & \frac{3|0|4}{0|3|4^{\mathstrut}} &
\frac{3|0|4|0}{5|0|6|0^{\mathstrut}} & \frac{3|4}{7|8^{\mathstrut}}\medskip\\
\frac{5|6|0}{0|1|2^{\mathstrut}} & \frac{5|6|0}{0|3|4^{\mathstrut}} &
\frac{5|6|0|0}{5|0|6|0^{\mathstrut}} & \frac{5|6}{7|8^{\mathstrut}}\medskip\\
\frac{7|8}{0|12^{\mathstrut}} & \frac{7|8}{0|34^{\mathstrut}} & \frac
{7|8|0|0}{5|0|6|0^{\mathstrut}} & \frac{7|8|0|0}{7|0|0|8^{\mathstrut}} 
\end{pmatrix}
.
\]
This matrix is almost identical to the matrix in Example \ref{ex22}, except
that the bipartition in the first column and second row is altered so that the
maximal equalizer is null and the BPM is indecomposable. The row dimension of
this matrix is $6$, the column dimension is $7$, and the entry dimension is
$23$. Therefore, the dimension is $36$.
\end{example}

\medskip

\section{Other results and future projects}

In Sections 4 we discussed the indecomposable factorization of a bipartition.
There is an inverse algorithm that recovers the bipartition from its
indecomposable factorization.

\begin{algorithm}
Let $c = \frac{B_{1} | B_{2} | \cdots| B_{r}}{A_{1} | A_{2} | \cdots|
A^{\mathstrut}_{r}}$ be a bipartition and let $c = C_{1} \cdots C_{r}$ be its
indecomposable factorization.

For $k$ from $1$ to $r:$

\indent\indent Let $n$ be the number of columns in $C_{k}$ and let $C_{k}
^{ij}$ be the $(i,j)^{th}$ entry of $C_{k}$.

\indent \indent Then $B_{k} = C^{11}_{k}$ and $A_{k} = C^{1n}_{k}$.
\end{algorithm}

The efficiency of the algorithm for computing the dimension of a BPM in
Appendix A can be improved by applying the \textquotedblleft
transpose-rotation.\textquotedblright

\begin{definition}
Let $c=\frac{B_{1}|B_{2}|\cdots|B_{r}}{A_{1}|A_{2}|\cdots|A_{r}^{\mathstrut}}$
be a bipartition. The \textbf{rotation of} $c$ is the bipartition
$c^{rot}:=\frac{A_{r}|A_{r-1}|\cdots|A_{1}}{B_{r}|B_{r-1}|\cdots|B_{1}}$. The
\textbf{transpose-rotation} of a BPM $C=(c_{ij})$ is the matrix $C^{T-R} 
:=(c_{ji}^{rot})$.
\end{definition}

The row and column dimensions of $C^{T-R}$ are the column and row dimensions
of $C$, respectively. Thus we can use the same algorithm to compute the row
and column dimensions of $C$. This computationally more efficient algorithm
appears in Appendix A.

Finally, there are \textquotedblleft generalized bipartition
matrices,\textquotedblright\ whose entries are bipartitions, products of
bipartition matrices, or products of generalized bipartition matrices. The
dimension algorithm for generalized bipartition matrices is identical to the
algorithm for bipartition matrices, but a computer implementation of the
algorithm in the more general setting has yet to be developed.

\section{Appendix A: The Python code}

\begin{lstlisting}[breaklines]
	import itertools
	
	def intersection(list1, list2):
	list3 = [v for v in list1 if v in list2]
	return list3
	
	def subsets(set):
	powerSet = []
	for i in range(1 << len(set)):
	powerSet.append([set[j] for j in range(len(set)) if (i & (1 << j))])
	return powerSet
	
	def setOfPowersets(r):
	p = []
	for s in r:
	p.append(subsets([i for i in range(1, len(s[1]))]))
	return p
	
	# returns all possible equalizers
	
	def possibleEqualizers(p):
	return list(itertools.product(*p))
	
	# checks if a matrix is equalized
	
	def equalized(e):
	for x in e:
	if x != e[0]:
	return 0
	return 1
	
	def combineEqualizers(e, isRow):
	if isRow == 1:
	return [transpose(transpose(eq)) for eq in list(itertools.product(*e))]
	else:
	return [transpose(eq) for eq in list(itertools.product(*e))]
	
	def transpose(m):
	return [[row[i] for row in m] for i in range(len(m[0]))]
	
	def unionize(m, i, j):
	unionizedList = []
	for h in range(i, j):
	unionizedList = unionizedList + m[h]
	return unionizedList
	
	# row equalizer for matrix of arbitrary size
	
	def arbitraryRE(b):
	equalizers = [[] for row in b]
	j = 0
	for row in b:
	pows = setOfPowersets(row)
	peqr = possibleEqualizers(pows)
	for eqr in peqr:
	i = 0
	eq = []
	for bip in row:
	n = 0
	eq.append([])
	for m in eqr[i]:
	if m < len(bip[1]) + 1:
	eq[i].append(unionize(bip[1], n, m))
	n = m
	elif n != len(bip[1]):
	eq[i].append(unionize(bip[1], n, len(bip[1])))
	for a in eq[i]:
	a.sort()
	i = i + 1
	if equalized(eq):
	equalizers[j].append(eqr)
	j = j + 1
	return equalizers
	
	# column equalizer for matrix of arbitrary size
	
	def arbitraryCE(b):
	bt = transpose(b)
	equalizers = [[] for col in bt]
	j = 0
	for col in bt:
	pows = setOfPowersets(col)
	peqr = possibleEqualizers(pows)
	for eqr in peqr:
	i = 0
	eq = []
	for bip in col:
	n = 0
	eq.append([])
	for m in eqr[i]:
	if m < len(bip[0]) + 1:
	eq[i].append(unionize(bip[0], n, m))
	n = m
	elif n != len(bip[0]):
	eq[i].append(unionize(bip[0], n, len(bip[0])))
	for a in eq[i]:
	a.sort()
	i = i + 1
	if equalized(eq):
	equalizers[j].append(eqr)
	j = j + 1
	return equalizers
	
	# returns all actual equalizers for a matrix
	
	def equalize(b):
	equalizers = []
	rowEqualizers = combineEqualizers(arbitraryRE(b), 1)
	colEqualizers = combineEqualizers(arbitraryCE(b), 0)
	for req in rowEqualizers:
	for ceq in colEqualizers:
	if req == ceq:
	equalizers.append(req)
	return equalizers
	
	# returns the maximal equalizer for a matrix
	
	def getMaximalEqualizer(b):
	e = equalize(b)
	maxEq = e[0].copy()
	for i in range(len(e)):
	for j in range(len(e[i])):
	jBreak = 0
	for k in range(len(e[i][j])):
	kBreak = 0
	if len(e[i][j][k]) > len(maxEq[j][k]):
	maxEq = e[i]
	jBreak = 1
	break
	elif len(e[i][j][k]) == len(maxEq[j][k]):
	for l in range(len(e[i][j][k])):
	if e[i][j][k][l] < maxEq[j][k][l]:
	maxEq = e[i]
	jBreak = 1
	kBreak = 1
	break
	else:
	break
	else:
	jBreak = 1
	kBreak = 1
	if kBreak == 1:
	break
	if jBreak == 0:
	break
	return maxEq
	
	# returns the embedding partition of a set and a subset
	
	def embeddingPartition(B, A):
	ACP = [[]]
	j = 0
	for i in range(len(B)):
	if B[i] in A:
	ACP[j].append(B[i])
	else:
	ACP.append([])
	j = j + 1
	return ACP
	
	# returns the indecomposable factorization of a bipartition
	
	def indecomposableFactorization(b):
	if len(b[0]) == 1:
	return [[[b]]]
	c = [[] for i in range(len(b[0]))]
	for k in range(len(b[0])):
	A = (embeddingPartition(unionize(b[0], 0, k + 1), b[0][k]))
	B = (embeddingPartition(unionize(b[1], k, len(b[0])), b[1][k]))
	for j in range(len(B)):
	c[k].append([])
	for i in range(len(A)):
	c[k][j].append([[A[i]], [B[j]]])
	return c
	
	# returns the indecomposable factorization of a BPM
	
	def factorization(b):
	if getMaximalEqualizer(b)[0][0] == []:
	return [b]
	factoredMatrix = [[[] for col in b[0]] for row in b]
	for row in range(len(b)):
	for col in range(len(b[row])):
	lam = getMaximalEqualizer(b)[row][col][0]
	A_ACP = []
	B_ACP = []
	Ai = []
	Bj = []
	setA = unionize(b[row][col][0], 0, len(b[row][col][0]))
	setA.sort()
	setB = unionize(b[row][col][1], 0, len(b[row][col][1]))
	setB.sort()
	A_ACP = (embeddingPartition(setA, unionize(b[row][col][0], lam, len(b[row][col][0]))))
	B_ACP = (embeddingPartition(setB, unionize(b[row][col][1], 0, lam)))
	for i in range(len(A_ACP)):
	Ai.append([intersection(A_ACP[i], b[row][col][0][l]) for l in range(lam, len(b[row][col][0]))])
	for i in range(len(B_ACP)):
	Bj.append([intersection(B_ACP[i], b[row][col][1][l]) for l in range(lam)])
	factoredMatrix[row][col] = [[[[[b[row][col][0][k] for k in range(lam)], [Bj[j][l] for l in range(lam)]]] for j in range(len(B_ACP))],
	[[[[Ai[i][k] for k in range(len(Ai[i]))], [b[row][col][1][l] for l in range(lam, len(b[row][col][1]))]] for i in range(len(A_ACP))]]]
	
	m1 = [[[[] for innerRow in factoredMatrix[row][0][0]], [[] for innerRow in factoredMatrix[row][0][1]]] for row in range(len(factoredMatrix))]
	for row in range(len(b)):
	for innerRow in range(len(factoredMatrix[row][0][0])):
	for col in range(len(b[row])):
	m1[row][0][innerRow] = m1[row][0][innerRow] + factoredMatrix[row][col][0][innerRow]
	for innerRow in range(len(factoredMatrix[row][0][1])):
	for col in range(len(b[row])):
	m1[row][1][innerRow] = m1[row][1][innerRow] + factoredMatrix[row][col][1][innerRow]
	
	m2 = [[], []]
	for row in range(len(b)):
	for innerRow in range(len(m1[row][0])):
	m2[0] = m2[0] + [m1[row][0][innerRow]]
	for innerRow in range(len(m1[row][1])):
	m2[1] = m2[1] + [m1[row][1][innerRow]]
	return [m2[0]] + factorization(m2[1])
	
	# calculates the row dimension of a BPM
	
	def rowDimension(b):
	dimension = 0
	for m in factorization(b):
	dimension = dimension + rowDimension2(m)
	return dimension
	
	def rowDimension2(b):
	rowDim = 0
	for row in b:
	isElementary = 1
	nonElementary = 1
	allEmpty = 1
	for bipartition in row:
	if len(bipartition[1]) > 1:
	isElementary = 0
	if len(bipartition[1]) <= 1:
	nonElementary = 0
	if nonElementary == 1:
	if getMaximalEqualizer([row]) != [[[] for bip in row]]:
	rowDim = rowDim + rowDimension([row])
	if isElementary == 0:
	for bipartition in row:
	if getMaximalEqualizer([row]) == [[[] for bip in row]]:
	for m in indecomposableFactorization(bipartition):
	rowDim = rowDim + rowDimension(m)
	if isElementary == 1:
	if len(row) >= 1:
	for bipartition in row:
	if bipartition[1] != [[]]:
	allEmpty = 0
	break
	if allEmpty == 1:
	isC = 0
	for bipartition in row:
	for s in bipartition[0]:
	isC = isC + len(s)
	if isC != 0:
	rowDim = rowDim + isC - 1
	return rowDim
	
	# calculates the column dimension of a BPM
	
	def colDimension(b):
	dimension = 0
	for m in factorization(b):
	dimension = dimension + colDimension2(m)
	return dimension
	
	def colDimension2(b):
	colDim = 0
	bt = transpose(b)
	for col in bt:
	nonElementary = 1
	isElementary = 1
	allEmpty = 1
	for bipartition in col:
	if len(bipartition[1]) > 1:
	isElementary = 0
	if len(bipartition[1]) <= 1:
	nonElementary = 0
	if nonElementary == 1:
	if getMaximalEqualizer(transpose([col])) != [[[]] for bip in col]:
	colDim = colDim + colDimension(transpose([col]))
	for bipartition in col:
	if len(bipartition[0]) > 1:
	isElementary = 0
	if getMaximalEqualizer(transpose([col])) == [[[]] for bip in col]:
	for m in indecomposableFactorization(bipartition):
	colDim = colDim + colDimension(m)
	if isElementary == 1:
	if len(col) >= 1:
	for bipartition in col:
	if bipartition[0] != [[]] and bipartition[0] != []:
	allEmpty = 0
	break
	if allEmpty == 1:
	osC = 0
	for bipartition in col:
	for s in bipartition[1]:
	osC = osC + len(s)
	if osC != 0:
	colDim = colDim + osC - 1
	return colDim
	
	# calculates the entry dimension of a BPM
	
	def entDimension(b):
	dimension = 0
	for m in factorization(b):
	dimension = dimension + entDimension2(m)
	return dimension
	
	def entDimension2(b):
	entDim = 0
	for row in b:
	for bipartition in row:
	if len(bipartition) == 2:
	if len(bipartition[0]) > 1:
	for m in indecomposableFactorization(bipartition):
	entDim = entDim + entDimension(m)
	elif len(bipartition[0]) == 1 and len(bipartition[1]) == 1 and len(bipartition[0][0]) != 0 and len(bipartition[1][0]) != 0:
	entDim = entDim + len(bipartition[0][0]) + len(bipartition[1][0]) - 1
	return entDim
	
	# calculates the dimension of a BPM
	
	def dimension(b):
	
	return rowDimension(b) + colDimension(b) + entDimension(b)
	
	# returns transpose-rotation of a matrix for the simplified dimension algorithm
	
	def newAlgorithm(b):
	bt_rot = []
	bt = transpose(b)
	for row in range(len(bt)):
	bt_rot.append([])
	for bipartition in range(len(bt[row])):
	bt_rot[row].append([[], []])
	for i in range(len(bt[row][bipartition][0])):
	bt_rot[row][bipartition][0].append(bt[row][bipartition][1][len(bt[row][bipartition][0]) - i - 1])
	bt_rot[row][bipartition][1].append(bt[row][bipartition][0][len(bt[row][bipartition][0]) - i - 1])
	return bt_rot
	
\end{lstlisting}

\section{Appendix B: Sample output}

\begin{lstlisting}[breaklines]
	print(factorization([[[[[], [], [], [1], [2]], [[1], [2], [3], [4], [5]]], [[[], [], [3], [4]], [[1], [2], [3, 4], [5]]], [[[], [], [5], [6]], [[1, 2], [4], [3], [5]]]]]))
	returns
	[[[[[[], []], [[1], [2]]], [[[], []], [[1], [2]]], [[[]], [[1, 2]]]], [[[[], []], [[], []]], [[[], []], [[], []]], [[[]], [[]]]], [[[[], []], [[], []]], [[[], []], [[], []]], [[[]], [[]]]], [[[[], []], [[], []]], [[[], []], [[], []]], [[[]], [[]]]]], [[[[[], [1]], [[3], [4]]], [[[3]], [[3, 4]]], [[[], [5]], [[4], [3]]]], [[[[], [1]], [[], []]], [[[3]], [[]]], [[[], [5]], [[], []]]]], [[[[[]], [[5]]], [[[2]], [[5]]], [[[]], [[5]]], [[[4]], [[5]]], [[[]], [[5]]], [[[6]], [[5]]]]]]
	
	print(dimension([[[[[], [], [], [1], [2]], [[1], [2], [3], [4], [5]]], [[[], [], [3], [4]], [[1], [2], [3, 4], [5]]], [[[], [], [5], [6]], [[1, 2], [4], [3], [5]]]]]))
	returns
	8
\end{lstlisting}
\medskip

\noindent \textbf{Acknowledgment.} We wish to thank Jim Stasheff for his helpful editorial suggestions and comments.

\end{document}